\documentclass[leqno]{article}
\usepackage{amsmath,amsthm}
\usepackage[psamsfonts]{amssymb}

\newtheorem{theorem}{Theorem}[section]
\newtheorem{proposition}[theorem]{Proposition}
\newtheorem{lemma}[theorem]{Lemma}

\theoremstyle{remark}
\newtheorem{remark}[theorem]{Remark}

\begin{document}

\title{\textbf{Natural connections on the bundle}\\
\textbf{of Riemannian metrics}}
\date{}
\author{\textsc{R. Ferreiro P\'erez
and
J. Mu\~noz Masqu\'e}\\
Insituto de F\'{\i}sica Aplicada, CSIC\\
C/ Serrano 144, 28006-Madrid\\
\emph{E-mail:}
\texttt{roberto@iec.csic.es, jaime@iec.csic.es}}

\maketitle

\begin{abstract}
\noindent Let $FM,\mathcal{M}_M$ be the bundles
of linear frames and Riemannian metrics
of a manifold $M$, respectively. The existence
of a unique $\mathrm{Diff}M$-invariant connection
form on
$J^1\mathcal{M}_M\times _MFM\to J^1\mathcal{M}_M$,
which is Riemannian with respect to the universal
metric on $J^1\mathcal{M}_M\times _MTM$, is proved.
Aplications to the construction of universal
Pontryagin and Euler forms, are given.
\end{abstract}

\bigskip


\smallskip


\smallskip

\noindent
\emph{Acknowledgements.\/} Supported by the Ministerio
de Educaci\'on y Ciencia of Spain, under grant
\#BFM2002--00141.

\section{Introduction}

Let $q\colon \mathcal{M}_M\to M$ be the bundle of
Riemannian metrics of a smooth manifold $M$ of dimension
$n$. The goal of this paper is to prove that the bundle
$q_1^\ast FM\to J^1\mathcal{M}_M$, obtained by
pulling the linear frame bundle $FM$ back to the
$1$-jet bundle of metrics, is endowed with a unique
$\mathrm{Diff}M$-invariant connection form
$\mathbf{\omega }$---called the universal Levi-Civita
connection---with the property of being Riemannian
with respect to the universal metric $\mathbf{g}$
on $q_1^\ast TM$; or equivalently, $\mathbf{\omega }$
is the only $\mathrm{Diff}M$-invariant connection form
on the subbundle $OM$ of pairs
$(u_x,j_x^1g)\in q_1^\ast FM$ such that $u_x$ is
$g_x$-orthonormal (see Theorem \ref{caracterizacion}
and the precise definitions below). This result is
analogous to that proving the existence of a canonical
connection on the principal $G$-bundle $J^1P\to C(P)$,
where $C(P)$ is the bundle of connections of a principal
$G$-bundle $P\to M$ (see \cite{geocon}).

As is well known (e.g., see \cite{Epstein}), the Levi-Civita
map, which assigns its Levi-Civita connection to every
Riemannian metric, is a natural map, i.e., it is
$\mathrm{Diff}M$-equivariant with respect to
the natural actions of the diffeomorphism group
on the space of Riemannian metrics on $M$ and on the space
of linear connections on $M$. This map induces a
$\mathrm{Diff}M$-invariant connection form
$\mathbf{\omega }_{\mathrm{hor}}$ on
$q_1^\ast FM\to J^1\mathcal{M}_M$, called the horizontal
Levi-Civita connection as it is horizontal with respect
to the projection
$q_1^\ast FM=J^1\mathcal{M}_M\times _MFM\to FM$;
but $\mathbf{\omega }_{\mathrm{hor}}$ is not a Riemannian
connection, i.e., it is not reducible to $OM$. Surprisingly,
$\mathbf{\omega }$ is obtained by adding a contact form to
$\mathbf{\omega }_{\mathrm{hor}}$, thus showing that the
contact structure on $J^1\mathcal{M}_M$ plays a crucial
role in our construction.

The connection form $\mathbf{\omega }$ allows us to construct
the universal Pontryagin and Euler differential forms on
$J^1\mathcal{M}_M$, which contain more information than
the corresponding cohomology classes on $M$; for example
the forms of degree greater than $n$ do not vanish
necessarily---unlike their cohomology classes. We also
remark on the fact that such forms play the same role,
in metric theory, than the universal characteristic forms
introduced in \cite{conn} in gauge theories.

\section{The geometry of the bundle of metrics}

\subsection{The bundle of metrics}

The bundle of Riemannian metrics
$q\colon \mathcal{M}_M\to M$ is a convex open subset in
$S^2(T^\ast M)$ and every Riemannian metric $g$ is
identified to a global section $g\colon M\to \mathcal{M}_M$
of this bundle.

Every system of coordinates $(U;x^i)$ on $M$ induces
a system of coordinates $(q^{-1}U;x^i,y_{ij})$ on
$\mathcal{M}_M$ by setting
$g_x=y_{ij}(g_x)(dx^i)_x\otimes (dx^j)_x$,
$\forall g_x\in \mathcal{M}_M$, $x\in U$. We denote by
$(y^{ij})$ the inverse matrix of $(y_{ij})$.

The diffeomorphism group of $M$ acts in a natural way
on $\mathcal{M}_M$ by automorphisms of this bundle:
The natural lift of a diffeomorphism
$\phi \in \mathrm{Diff}M$ to the bundle of metrics
$\bar{\phi }\colon \mathcal{M}_M\to \mathcal{M}_M$
is defined by
\begin{equation}
\bar{\phi }(g_x)
=\left( \phi ^\ast
\right) ^{-1}(g_x)
\in (\mathcal{M}_M)_{\phi (x)},
\label{barphi}
\end{equation}
$\phi ^\ast \colon S^2T_{\phi (x)}^\ast M
\to S^2T_x^\ast M$ being the induced homomorphism.
Hence $q\circ \bar{\phi }=\phi \circ q$.
In the same way, the lift of a vector field
$X\in \frak{X}(M)$ is denoted by
$\bar{X}\in \frak{X}(\mathcal{M}_M)$. If
$X=X^i\partial /\partial x^i$, then
\begin{equation}
\bar{X}
=X^i
\frac{\partial }{\partial x^i}
-\sum_{i\leq j}
\left(
\frac{\partial X^r}{\partial x^i}
y_{kj}
+\frac{\partial X^k}{\partial x^j}
y_{ki}
\right)
\frac{\partial }{\partial y_{ij}}.
\label{expXbarra}
\end{equation}

\subsection{Jets of metrics}

Let $q_r\colon J^r\mathcal{M}_M\to M$ be the $r$-jet bundle
of sections of $\mathcal{M}_M$ and, for every $r\geq s$, let
$q_{rs}\colon J^r\mathcal{M}_M\to J^s\mathcal{M}_M$ be
the canonical projections. For every $\phi \in \mathrm{Diff}M$
we denote by $\bar{\phi }^{(r)}$ the natural prolongation to
$J^r\mathcal{M}_M$ of the lift $\bar{\phi }$ given in
\eqref{barphi}; precisely,
$\bar{\phi }^{(r)}(j_x^rg)=j_x^r
\left(
 \bar{\phi }\circ g\circ \phi^{-1}
\right) $.
Similarly, $\bar{X}^{(r)}$ denotes the jet prolongation
 of the lift $\bar{X}\in\frak{X}(\mathcal{M}_M)$ given in
\eqref{expXbarra}.

Let $(q_{r}^{-1}U;x^i,y_{ij,}y_{ij,I})$, $1\leq|I|\leq r$,
$I\in \mathbb{N}^n$, be the coordinate system induced by
$(q^{-1}U;x^i,y_{ij})$; i.e., $y_{ij,I}(j_x^1g)
=(\partial^{|I|}(y_{ij}\circ g)/\partial x^I)(x)$.
If $(U;x^i)$ is a normal coordinate system for the metric
$g$ centered at $x$, then we have
$y_{ij}(j_x^1g)=\delta _{ij}$, $y_{ij,k}(j_{x}^1g)=0$.

Let us fix a coordinate system $(U;x^i)$ centered at $x\in M$,
and let $\phi \in\mathrm{Diff}M$ be a diffeomorphism such that
$\phi(x)\in U$. The equations of the transformation
$\bar{\phi }^{(1)}$ are as follows:

\begin{equation}
\left\{
\begin{array}
[c]{l}
(\phi ^{-1})^a=x^a\circ \phi^{-1},\\
y_{ij}\circ \bar{\phi}^{(1)}
=y_{ab}
\left(
\dfrac{\partial(\phi^{-1})^a}{\partial x^i}
\circ\phi
\right)
\left(
\dfrac{\partial(\phi^{-1})^b}{\partial x^j}
\circ\phi
\right) ,\\
\\
y_{ij,k}\circ \bar{\phi }^{(1)}
=y_{ab,c}
\left(
\dfrac{\partial (\phi^{-1})^{c}}{\partial x^k}
\circ\phi
\right)
\left(
\dfrac{\partial (\phi^{-1})^a}{\partial x^i}
\circ \phi
\right)
\left(
\dfrac{\partial (\phi^{-1})^b}{\partial x^j}
\circ\phi
\right) \\
\qquad \qquad \qquad
+y_{ab}
\left(
\dfrac{\partial (\phi^{-1})^b}{\partial x^j}
\circ \phi
\right)
\left(
\dfrac{\partial ^2(\phi ^{-1})^a}
{\partial x^i\partial x^k}
\circ \phi
\right) \\
\qquad \qquad \qquad
+y_{ab}
\left(
\dfrac{\partial (\phi^{-1})^a}{\partial x^i}
\circ\phi
\right)
\left(
\dfrac{\partial ^2(\phi^{-1})^b}
{\partial x^j\partial x^k}
\circ \phi
\right) .
\end{array}
\right. \label{coor}
\end{equation}

If $M$ is an orientable and connected manifold,
we denote by $\mathrm{Diff}^{+}M$ the subgroup
of orientation preserving diffeomorphisms.
The manifold is
said to be irreversible if
$\mathrm{Diff}M=\mathrm{Diff}^{+}M$;
otherwise, the manifold is said to be reversible.
Every compact, orientable and connected
manifold of dimension $\leq 3$ is reversible; e.g.,
see \cite[Chapter 9, \S 1]{Hirsch}. We recall
the following result about extending diffeomorphisms
(see \cite{Hirsch, Palais}):

\begin{lemma}
\label{dext}Let $\phi\colon U\to U$ be a diffeomorphism
defined on an open neighbourhood of $x$ in an orientable
differentiable manifold $M$ such that $\phi(x)=x$. If $M$
is irreversible we further assume
$\phi _{\ast ,x}\in Gl^{+}(T_xM)$. Then a global diffeomorphism
$\tilde{\phi }\in \mathrm{Diff}M$ exists coinciding with
$\phi $ on a neighbourhood of $x$.
\end{lemma}

$\mathrm{Diff}^{+}M$ (and hence $\mathrm{Diff}M$) acts
transitively on every orientable connected manifold $M$.
Even more, as a simple consequence of the existence of normal
coordinates and Lemma \ref{dext} we have

\begin{proposition}
\label{trans}If $M$ is an orientable and connected
manifold, then the group $\mathrm{Diff}^{+}M$ acts
transitively on $J^1\mathcal{M}_M$.
\end{proposition}

Let $\theta \in \Omega ^1
(J^1\mathcal{M}_M,S^2(T^\ast M))$ be the
structure form of $J^1\mathcal{M}_M$ (see \cite{Mu}),
 where we use the canonical identification
$V(\mathcal{M}_M)\simeq \mathcal{M}_M
\times _MS^2(T^\ast M)$; i.e.,
$\theta _{j_x^1g}(X)=(q_{10})_\ast (X)
-g_\ast (q_{1\ast}(X))$,
$\forall X\in T_{j_x^1g}J^1\mathcal{M}_M$.
In local coordinates,
\begin{equation}
\theta
=(dy_{ij}-y_{ij,k}dx^k)\otimes dx^i\otimes dx^j.
\label{exptheta}
\end{equation}

The bundle
$q_1^\ast TM=J^1\mathcal{M}_M\times _MTM
\to J^1\mathcal{M}_M$, obtained by pulling $TM$ back
via $q_1\colon J^1\mathcal{M}_M\to M$, is endowed
with a universal metric given by
$\mathbf{g}
\left(
\left(
j_{x}^1g,X
\right) ,
\left(
j_{x}^1g,Y
\right)  \right)
=g_{x}(X,Y)$, $\forall X,Y\in T_{x}M$, which
satisfies the following universal property:
$(j^1g)^\ast \mathbf{g}=g$, for every Riemannian metric
$g$ on $M$. By means of this metric, we can identify
$q_1^\ast TM$ to $q_1^\ast T^\ast M$.

If $\alpha \in\Omega ^\bullet
(J^1\mathcal{M}_M,T^\ast M)$, the element of
$\Omega ^\bullet (J^1\mathcal{M}_M,TM)$
corresponding to $\alpha $ under this identification,
is denoted by $\mathbf{g}^{-1}\alpha$. We have
\[
q_1^\ast
\left(
\mathrm{End}TM
\right)
\simeq q_1^\ast
\left(
\otimes ^2T^\ast M
\right)
=q_1^\ast
\left(
S^{2}T^\ast M
\right)
\oplus q_1^\ast
\left(
\wedge ^2T^\ast M
\right) .
\]
Let $\mathrm{End}_{\text{\textsc{S}}}TM$
(resp.\
$\mathrm{End}_{\text{\textsc{A}}}TM$)
be the image of $q_1^\ast
\left(
S^{2}T^\ast M
\right) $
(resp.\ $q_1^\ast
\left(
\wedge ^2T^\ast M
\right) $)
in $q_1^\ast
\left(  \mathrm{End}TM
\right)$
under the previous isomorphism.
We have
\begin{equation}
\vartheta
=\mathbf{g}^{-1}\theta
\in \Omega ^1(J^1\mathcal{M}_M,
\mathrm{End}_{\text{\textsc{S}}}TM).
\label{gtheta}
\end{equation}
If
$\alpha
\in \Omega ^r(J^1\mathcal{M}_M,\mathrm{End}TM)$,
then a decomposition exists such that
$\alpha
=\alpha _{\text{\textsc{S}}}
+\alpha _{\text{\textsc{A}}}$,
where the forms
$\alpha _{\text{\textsc{S}}}\in \Omega ^r
(J^1\mathcal{M}_M,\mathrm{End}_{\text{\textsc{S}}}TM)$,
$\alpha _{\text{\textsc{A}}}\in \Omega ^r
(J^1\mathcal{M}_M,\mathrm{End}_{\text{\textsc{A}}}TM)$
are called the symmetric and anti-symmetric parts of
$\alpha $, respectively.

\section{Natural connections}

\subsection{Linear frame bundles}

Let $\pi \colon FM\to M$ be the linear frame bundle
of $M$, and let $(x^h,x_j^i)$ be the coordinate system
induced on $\pi^{-1}U$ by a coordinate system $(U;x^h)$
in $M$; i.e., $u=((\partial /\partial x^1)_x,\dotsc,
(\partial /\partial x^n)_x)\cdot(x_j^i(u))$,
$\forall u\in\pi ^{-1}(x)$, $\forall x\in U$.
The lift of $\phi \in \mathrm{Diff}M$ to $FM$
is denoted by $\tilde{\phi }\colon FM\to FM$,
$\tilde{\phi }(u)=\phi _\ast (u)$. Analogously,
$\tilde{X}\in \frak{X}(FM)$ stands for the lift
of $X\in\frak{X}(M)$. Let
$q_1^\ast FM=J^1\mathcal{M}_M\times _MFM$
be the pull-back of $FM$ to $J^1\mathcal{M}_M$
via $q_1\colon J^1\mathcal{M}_M\to M$.
There are two canonical projections
\[
\begin{array}
[c]{ccc}
q_1^\ast FM
& \overset{\bar{q}_1}{\longrightarrow }
& FM\\
{\scriptstyle \bar{\pi }}\downarrow
&  &
\downarrow{ \scriptstyle\pi }\\
J^1\mathcal{M}_M
&
\overset{q_1}{\longrightarrow }
&
M
\end{array}
\]
The first projection
${\bar{\pi }}\colon q_1^\ast FM\to J^1\mathcal{M}_M$
is a principal $Gl(n,\mathbb{R})$-bundle with respect
to the induced action, given by
$(j_x^1g,u)\cdot A=(j_x^1g,u\cdot A)$,
$\forall j_x^1g\in J_x^1\mathcal{M}_M$,
$\forall u\in F_xM$,
$\forall A\in Gl(n,\mathbb{R})$, and the second
projection $\bar{q}_1\colon q_1^\ast FM\to FM$
is $Gl(n,\mathbb{R})$-equivariant.

The diffeomorphism group $\mathrm{Diff}M$ acts on
$\mathcal{M}_M$ and on $FM$ as explained above;
hence it acts on $q_1^\ast FM$ by the induced
action. If $\phi\in\mathrm{Diff}M$ its lift to
$q_1^\ast FM$ is $\hat{\phi }
=(\bar{\phi }^{(1)},\tilde{\phi })$.
Similarly, if $X\in \frak{X}(M)$ we denote by
$\hat{X}$ its lift to $q_1^\ast FM$. As $\bar{q}_1$
is a $\mathrm{Diff}M$-equivariant map we have
$(\bar{q}_1)_\ast (\hat{X})=\tilde{X}$.
For every $t\in\mathbb{R}$, we define
$\varphi _t\in \mathrm{Aut}FM$ (resp.\
$\bar{\varphi }_t
\in \mathrm{Diff}(J^1\mathcal{M}_M)$,
resp.\ $\hat{\varphi }_t\in \mathrm{Aut}
\left(
q_1^\ast FM
\right) $)
by
\begin{align*}
\varphi _t(u)
&
=\exp( -\tfrac{t}{2})\cdot u,\\
\bar{\varphi }_t(j_x^1g)
&
=j_x^1
\left(
\exp (t)\cdot g
\right) ,\\
\hat{\varphi }_t(j_x^1g,u)
&
=\left(
j_x^1
\left(
\exp (t)\cdot
g\right) ,
\exp (-\tfrac{t}{2})\cdot u
\right) .
\end{align*}
We denote by $\xi \in \frak{X}(FM)$
(resp.\ $\mathbf{\xi }\in \frak{X}(J^1\mathcal{M}_M)$,
resp.\ $\mathbf{\hat{\xi }}\in \frak{X}(q_1^\ast FM)$)
the infinitesimal generator of the $1$-parameter group
$(\varphi _t)$ (resp.\ $(\bar{\varphi }_t)$,
resp.\ $(\hat{\varphi }_t)$) defined above. We have
$q_{1\ast }(\mathbf{\xi })=0$,
$\bar{q}_{1\ast }(\mathbf{\hat{\xi })}=\xi $, and
$\bar{\pi }_\ast (\mathbf{\hat{\xi })}=\mathbf{\xi }$.
Hence, the group
$\mathcal{G}=\mathrm{Diff}M\times\mathbb{R}$
acts on the principal $Gl(n,\mathbb{R})$-bundle
$\bar{\pi }\colon q_1^\ast FM\to J^1\mathcal{M}_M$.
by automorphisms. If $(\phi ,t)\in \mathcal{G}$,
we set $\hat{\phi }_t
=\hat{\phi }\circ \hat{\varphi }_t
=\hat{\varphi }_t\circ \hat{\phi }$.
This action induces a $\mathcal{G}$-action
on the associated bundles to $q_1^\ast FM$, (such as
$q_1^\ast TM$, $q_1^\ast T^\ast M$, etc.), and
on the space of sections and differential forms with
values on such bundles.

\begin{proposition}
\label{ginv}
The universal metric $\mathbf{g}\in \Omega ^0
(J^1\mathcal{M}_M,S^2T^\ast M)$ on the bundle
$q_1^\ast TM$, is invariant under the action
of the group
$\mathcal{G}=\mathrm{Diff}M\times\mathbb{R}$
defined above.
\end{proposition}

\begin{proof}
Let
$\mathbf{\bar{g}}\in \Omega ^0(q_1^\ast FM,S^2
\left(
\mathbb{R}^n
\right) ^\ast )$
be the $Gl(n,\mathbb{R})$-invariant function on
$q_1^\ast FM$ corresponding to $\mathbf{g}$.
Let $(e_i)$ be the standard basis on
$\mathbb{R}^n$ and let $(e^i)$ be its dual basis.
For every $j_x^1g\in J^1\mathcal{M}_M$
and every frame $u_x
=\left(
X_1,\dotsc,X_n
\right)
\in FM_x$ we have $\mathbf{\bar{g}}(j_x^1g,u_x)
=g_x(X_i,X_j)e^i\otimes e^j$. Hence,
for every $(\phi ,t)\in\mathcal{G}$, we have
\begin{align*}
\left(
\hat{\phi}_{t}^\ast \mathbf{\bar{g}}
\right)
(j_{x}^1g,u_{x})
&
=\mathbf{\bar{g}}
(\bar{\phi }_t^{(1)}
(j_x^1g),\tilde{\phi }_tu_x)\\
&
=\left(
\phi_{t}^{-1}
\right) ^\ast g_{\phi (x)}
(\phi _{t\ast }X_i,\phi _{t\ast }X_j)
e^i\otimes e^j\\
&
=g_x(X_i,X_j)e^i\otimes e^j\\
&
=\mathbf{\bar{g}}(j_x^1g,u_x).
\end{align*}
\end{proof}

\subsection{The horizontal Levi-Civita connection}

We denote by $\pi\colon F_{g}M\to M$ the orthonormal
frame bundle with respect to a Riemannian metric $g$ on
$M$, which is a reduction of group $O(n)$ of the bundle
$FM$. We denote by $\Gamma ^g$ the Levi-Civita
connection of $g$; i.e., the only symmetric connection on
$F_{g}M$. If there is no risk of confusion we also denote
by $\Gamma ^g$ its direct image with respect to
the canonical injection $F_{g}M\hookrightarrow FM$
(see \cite[II. Proposition 6.1]{KN}). Analogously,
$\omega ^g$ denotes the connection form of both connections,
and $\nabla ^g$ the covariant derivation law with
respect to $\Gamma ^g$ on the associated vector bundles.

\begin{proposition}
The $\frak{gl}(n,\mathbb{R})$-valued $1$-form on
$q_1^\ast FM$ defined by
\[
\mathbf{\omega }_{\mathrm{hor}}(X)
=\omega ^g((\bar{q}_1)_\ast X),\;
\forall X\in T_{(j_x^1g,u)}(q_1^\ast FM),
\]
is a connection form on the principal
$Gl(n,\mathbb{R})$-bundle ${\bar{\pi
}\colon}q_1^\ast FM\to J^1\mathcal{M}_M$.
\end{proposition}

\begin{proof}
The definition makes sense as
$\omega ^g|_{\pi ^{-1}(x)}$ depends only on
$j_x^1g$. We check the two characteristic
properties of a connection form.

(1) For every $A\in \frak{gl}(n,\mathbb{R})$
we have
\[
\begin{array}
[c]{lll}
\mathbf{\omega }_{\mathrm{hor}}
(A_{(j_x^1g,u)}^\ast )
&
\!
=\omega ^g
((\bar{q}_1)_\ast A_{(j_{x}^1g,u)}^\ast )
&
\text{by the definition of }
\mathbf{\omega }_{\mathrm{hor}}\\
&
\!
=\omega^{g}(A_{u}^\ast )
&
\text{as }\bar{q}_1
\text{ is equivariant}\\
&
\!
=A
&
\text{as }\omega ^g
\text{ is a connection form}
\end{array}
\]
(2) For every $A\in Gl(n,\mathbb{R})$ and
$X\in TFM_{J^1\mathcal{M}_M}$, we have

\[
\begin{array}
[c]{lll}
(R_A^\ast \mathbf{\omega }_{\mathrm{hor}})(X)
&
\!
=\mathbf{\omega }_{\mathrm{hor}}
((R_A)_\ast X)
& \\
&
\!
=\omega^{g}((\bar{q}_1)_\ast
(R_A)_\ast X)
& \\
&
\!
=\omega ^g((R_A)_\ast (\bar{q}_1)_\ast X)
&
\text{as }
\bar{q}_1
\text{ is equivariant}\\
&
\!
=(R_A^\ast \omega ^g)
((\bar{q}_1)_\ast X)
& \\
&
\!
=(\mathrm{Ad}_{A^{-1}}\circ \omega ^g)
((\bar{q}_1)_\ast X)
&
\text{as }\omega ^g
\text{ is a conection form}\\
&
\!
=(\mathrm{Ad}_{A^{-1}}
\circ \mathbf{\omega }_{\mathrm{hor}})(X).
&
\end{array}
\]
\end{proof}

\begin{remark}
The bundle of orthonormal frames cannot be used
in the preceding definition as it depends
on the metric chosen. Below, we show that,
in fact, $\mathbf{\omega }_{\mathrm{hor}}$
is not reducible to the bundle of orthonormal frames,
although a new connection form $\mathbf{\omega }$
can be defined, which will be reducible to this bundle.
\end{remark}

The connection
$\mathbf{\omega }_{\mathrm{hor}}$ induces a derivation
law $\nabla ^{\mathbf{\omega }_{\mathrm{hor}}}$
in the associated bundles to $q_1^\ast FM$;
in particular, on
$q_1^\ast TM$, $q_1^\ast T^\ast M$,
etc. In local coordinates we have
\begin{align*}
\nabla^{\mathbf{\omega }_{\mathrm{hor}}}X
&
=\left(
dX^i+
\text{\mbox{\boldmath$\Gamma $}}_{jk}^iX^kdx^j
\right)
\otimes \frac{\partial}{\partial x^i},\\
\nabla ^{\mathbf{\omega }_{\mathrm{hor}}}\alpha
&
=\left(
d\alpha_i-\text{\mbox{\boldmath$\Gamma $}
}_{ji}^k\alpha_kdx^j
\right)
\otimes dx^i,
\end{align*}
\begin{equation}
\mbox{\boldmath$\Omega $}_{\mathrm{hor}}
=\left(
d\mbox{\boldmath$\Gamma $}_{jk}^i\wedge dx^k
+\mbox{\boldmath$\Gamma $}_{as}^i
\mbox{\boldmath$\Gamma $}_{jr}^adx^s\wedge dx^r
\right)
dx^j\otimes
\frac{\partial}{\partial x^i},
\label{Omegahor}
\end{equation}
where
$X=X^i\partial /\partial x^i
\in \Omega ^0(J^1\mathcal{M}_M,TM)$,
$\alpha =\alpha _idx^i
\in \Omega ^0(J^1\mathcal{M}_M,T^\ast M)$,
\begin{equation}
\text{\mbox{\boldmath$\Gamma$}}_{jk}^i
=\tfrac{1}{2}y^{ia}(y_{aj,k}+y_{ak,j}-y_{jk,a}),
\label{Gammanegrita}
\end{equation}
and $\mbox{\boldmath$\Omega $}_{\mathrm{hor}}$
is the curvature form of
$\mathbf{\omega }_{\mathrm{hor}}$.

\begin{proposition}
\label{propomega}The connection form
$\mathbf{\omega }_{\mathrm{hor}}$
satisfies the following properties:

\begin{enumerate}
\item [\emph{(1)}]If
$\sigma_{g}\colon FM\to q_1^\ast FM$ is the
equivariant section induced by a Riemannian metric
$g$ (i.e., $\sigma _g(u_x)=(j_x^1g,u_x)$), then
$\sigma _g^\ast \mathbf{\omega }_{\mathrm{hor}}$
is the Levi-Civita connection form of $g$.

\item[\emph{(2)}] The form
$\mathbf{\omega }_{\mathrm{hor}}$ is invariant under
the action of the group
$\mathcal{G}=\mathrm{Diff}M\times\mathbb{R}$ on
$q_1^\ast FM$.

\item[\emph{(3)}]
$\nabla ^{\mathbf{\omega }_{\mathrm{hor}}}\mathbf{g}
=\theta $.
\end{enumerate}
\end{proposition}

\begin{proof}
(1) Let $X\in T_{u_x}FM$. As
$\bar{q}_1\circ\sigma _g=\mathrm{id}_{FM} $,
we have
\begin{align*}
(\sigma _g^\ast \mathbf{\omega }_{\mathrm{hor}})(X)
& =\mathbf{\omega }_{\mathrm{hor}}
(\sigma _{g\ast }(X))\\
& =\omega ^g(\bar{q}_{1\ast }\sigma _{g\ast }(X))\\
& =\omega ^g
\left(
\left(
\bar{q}_1\circ \sigma _g
\right) _\ast (X)
\right) \\
&
=\omega ^g(X).
\end{align*}

(2) First recall that if $\phi\in\mathrm{Diff}M$ and
$\omega^{g}$ is the Levi-Civita connection of the metric
$g$, then $(\tilde{\phi }^{-1})^\ast \omega ^g
=\omega ^{\phi \cdot g}$ is the Levi-Civita connection
of the metric $\phi \cdot g
=\bar{\phi }\circ g\circ \phi ^{-1}$.
In the same way, if $t\in\mathbb{R}$, then
$(\tilde{\varphi }_t^{-1})^\ast \omega ^g$ is the
Levi-Civita connection of the metric $\exp (t)\cdot g$.
From the definition of $\mathbf{\omega }_{\mathrm{hor}}$,
for every $X\in T_{(j_x^1g,u)}q_1^\ast FM$, and
$(\phi ,t)\in\mathcal{G}$ we have
\begin{align*}
(\hat{\phi }_t^\ast
\mathbf{\omega }_{\mathrm{hor}})(X)
&
=(\mathbf{\omega }_{\mathrm{hor}})
(\hat{\phi }_{t\ast }X)\\
&
=\omega ^{\phi\cdot g}
(\bar{q}_{1\ast}\hat{\phi }_{t\ast }X)\\
&
=\left(
\left(
\tilde{\phi }_t^{-1}
\right) ^\ast \omega ^g
\right)
(\tilde{\phi }_{t\ast}\bar{q}_{1\ast }X)\\
&
=\omega ^g(\bar{q}_{1\ast }X)\\
&
=\mathbf{\omega }_{\mathrm{hor}}(X).
\end{align*}

(3) In local coordinates, we have
\begin{align*}
\nabla ^{\mathbf{\omega }_{\mathrm{hor}}}
\mathbf{g}
&
=\left(
dy_{ij}
-\left(
y_{aj}^a\mbox{\boldmath$\Gamma $}_{ki}
+y_{ai}^a\mbox{\boldmath$\Gamma $}_{kj}
\right)
dx^k
\right)
\otimes dx^j\otimes dx^i\\
&
=\left(
dy_{ij}-y_{ij,k}dx^k
\right)
\otimes dx^j\otimes dx^i,
\end{align*}
where we have used the equation
\begin{equation}
y_{ij,k}
=y_{aj}\mbox{\boldmath$\Gamma $}_{ki}^a
+y_{ai}\mbox{\boldmath$\Gamma $}_{kj}^a,
\label{crist}
\end{equation}
and we conclude by virtue of the formula
\eqref{exptheta}.
\end{proof}

\subsection{The universal Levi-Civita connection}

The connections on $q_1^\ast FM$ are an affine space
modelled over
$\Omega ^1(J^1\mathcal{M}_M,$\textrm{End}$TM)$.
Furthermore, as
$\vartheta \in \Omega ^1
(J^1\mathcal{M}_M,\mathrm{End}TM)$,
we can define a connection form on $q_1^\ast FM$
as follows:
\[
\mathbf{\omega }=\mathbf{\omega }_{\mathrm{hor}}
+\tfrac{1}{2}\vartheta .
\]
The connection form $\mathbf{\omega }$ is called the
universal Levi-Civita connection.

The following lemma can be proved by computing
in local coordinates:

\begin{lemma}
\label{lemasym}
If
$\alpha \in \Omega ^1
(J^1\mathcal{M}_M,\mathrm{End}TM)$,
then
$\nabla^{\mathbf{\omega }_{\mathrm{hor}}
+\alpha}\mathbf{g}
=\nabla ^{\mathbf{\omega }_{\mathrm{hor}}}\mathbf{g}
-2\alpha_{\text{\textsc{S}}}$.
\end{lemma}

The next theorem states the basic properties
of the universal Levi-Civita connection
and it is analogous to Proposition \ref{propomega}.

\begin{theorem}
\label{propomeganegrita}The connection form
$\mathbf{\omega }$ satisfies the following
properties:

\begin{enumerate}
\item [\emph{(1)}]With the notations of
\emph{Proposition \ref{propomega}},
for any Riemannian metric $g$, the form
$\sigma _g^\ast \mathbf{\omega }$ is
the Levi-Civita connection form of $g$.

\item[\emph{(2)}] The form $\mathbf{\omega }$
is invariant under the action of
$\mathcal{G}=\mathrm{Diff}M\times \mathbb{R}$
on $q_1^\ast FM$.

\item[\emph{(3)}] If $\nabla ^{\mathbf{\omega }}$
is the derivation law induced by $\mathbf{\omega }$,
then
\begin{equation}
\nabla^{\mathbf{\omega }}\mathbf{g}=0.
\label{nablaG}
\end{equation}
\end{enumerate}
\end{theorem}

\begin{proof}
(1) Follows from Proposition \ref{propomega}
and from the fact that $(j^1g)^\ast \vartheta =0$.

(2) From Proposition \ref{ginv} and
Proposition \ref{propomega}--(2), we know
that $\mathbf{\omega }_{\mathrm{hor}}$ and
$\mathbf{g}$ are $\mathcal{G}
$-invariant. Moreover, from
Proposition \ref{propomega}--(3) we have
$\theta
=\nabla ^{\mathbf{\omega }_{\mathrm{hor}}}\mathbf{g}$,
and hence $\theta $, as well as
$\vartheta =\mathbf{g}^{-1}\theta $,
are also $\mathcal{G}$-invariant. Hence $\mathbf{\omega }
=\mathbf{\omega }_{\mathrm{hor}}
+\tfrac{1}{2}\vartheta$ is $\mathcal{G}$-invariant.

(3) It is a consequence of Lemma \ref{lemasym}
and Proposition \ref{propomega}--(3).
\end{proof}

Let $OM\to J^1\mathcal{M}_M$ be the reduction
of $q_1^\ast FM $ to the subgroup $O(n)$ given
by $OM=
\left\{
(u_{x},j_{x}^1g)\in q_1^\ast FM\colon u_x\,
\text{is }g_x\text{-orthonormal}
\right\} $.

The following result shows the advantage
of $\mathbf{\omega }$ over
$\mathbf{\omega }_{\mathrm{hor}}$.

\begin{proposition}
\label{reducible}
The connection $\mathbf{\omega }$ is reducible
to a connection on $OM$.
\end{proposition}

\begin{proof}
It follows from \cite[III, Proposition 1.5]{KN}
and the formula \eqref{nablaG}.
\end{proof}

\begin{proposition}
If $\mbox{\boldmath$\Omega $}$ is the curvature
form of $\mathbf{\omega }$, then we have
\begin{equation}
\mbox{\boldmath$\Omega $}
=\left(
\mbox{\boldmath$\Omega$}_{\mathrm{hor}}
\right) _{\text{\textsc{A}}}
-\tfrac{1}{2}\vartheta \wedge \vartheta .
\label{Omegan}
\end{equation}
\end{proposition}

\begin{proof}
As $\mathbf{\omega }
=\mathbf{\omega }_{\mathrm{hor}}
+\frac{1}{2}\vartheta $,
we have $\mbox{\boldmath$\Omega $}
=\mbox{\boldmath$\Omega $}_{\mathrm{hor}}
+\tfrac{1}{2}
d^{\mathbf{\omega }_{\mathrm{hor}}}\vartheta
+\tfrac{1}{4}\vartheta \wedge \vartheta $.
Hence, it suffices to prove
$d^{\mathbf{\omega }_{\mathrm{hor}}}\vartheta
=-\vartheta \wedge \vartheta
-2\left(
\mbox{\boldmath$\Omega $}_{\mathrm{hor}}
\right) _{\text{\textsc{S}}}$. Let $(x^i)$
be a system of
normal coordinate for $g$ centered at $x$.
By taking the covariant exterior
differential with respect to
$\mathbf{\omega }_{\mathrm{hor}}$ in the local
expression $\vartheta
=y^{ia}(dy_{aj}-y_{aj,k}\wedge dx^k)
\otimes dx^j\otimes \partial /\partial x^i$,
and evaluating it at $j_x^1g$, we have
\begin{equation}
\left(
d^{\mathbf{\omega }_{\mathrm{hor}}}\vartheta
\right) _{j_{x}^1g}
=\left(
dy^{ia}\wedge dy_{aj}-dy_{ij,k}\wedge dx^k
\right) _{j_x^1g}\otimes
\left(
dx^j\otimes
\frac{\partial }{\partial x^i}
\right) _{j_x^1g}.
\label{dertheta}
\end{equation}
Taking the exterior differential in
$y^{ia}y_{aj}=\delta_j^i$ we obtain
$dy^{ia}=-y^{ib}dy_{bj}y^{jb}$, and taking
the exterior differential in the formula \eqref{crist}
and evaluating at $j_x^1g$ we have
$\left(
dy_{ij,k}
\right) _{j_x^1g}
=(d\mbox{\boldmath$\Gamma $}_{jk}^i)_{j_x^1g}
+(d\mbox{\boldmath$\Gamma $}_{ik}^j)_{j_x^1g}$.
Substituting these expressions in \eqref{dertheta}
and taking the formula \eqref{Omegahor}
into account we obtain
\begin{align*}
\left(
d^{\mathbf{\omega }_{\mathrm{hor}}}\vartheta
\right) _{j_{x}^1g}\!
&
=
\!-\left(
dy_{ia}\!\wedge \!dy_{aj}\!
+\!
\left(
d\mbox{\boldmath$\Gamma $}_{jk}^i\!
+\! d\mbox{\boldmath$\Gamma $}_{ik}^j
\right)
\!\wedge\!dx^k
\right) _{j_x^1g}\!\otimes\!
\left(
\!dx^j\!\otimes\!
\frac{\partial}{\partial x^i}\!
\right) _{j_{x}^1g}\\
\!
&
=\!-\vartheta _{j_x^1g}\wedge \vartheta _{j_x^1g}
-2\left(
\left(
\mbox{\boldmath$\Omega $}_{\mathrm{hor}}
\right) _{\text{\textsc{S}}}
\right) _{j_{x}^1g}.
\end{align*}
\end{proof}

\section{Universal Pontryagin and Euler forms}

Let $\mathcal{I}_{d}^{G}$ denote the Weil invariant
polynomials of degree $d$ for the Lie group $G$,
see \cite[XII]{KN}. As
$OM\to J^1\mathcal{M}_M$ is a principal
$O(n)$-bundle and $\mathbf{\omega }$ is a connection
form on this bundle, the Chern-Weil construction of the
characteristic classes provides us a closed differential
$(2d)$-form $f(\mbox{\boldmath$\Omega $})$ on
$J^1\mathcal{M}_M$, by applying a Weil polynomial
$f\in \mathcal{I}_{d}^{O(n)}$ to the curvature
$\mbox{\boldmath$\Omega $}$ of $\mathbf{\omega }$.
As is well known (e.g., see \cite{GS, KN}),
$\mathcal{I}^{O(n)}$ is spanned by the polynomials
$p_k\in\mathcal{I}_{2k}^{O(n)}$,
$1\leq k\leq [\frac{n}{2}]$, characterized by
\[
\det
\left(
\lambda I-\tfrac{1}{2\pi}X
\right)
=\sum p_k(X)\lambda ^{n-2k},
\quad \forall X\in \frak{so}(n).
\]
We define the universal $k$-Pontryagin form of $M$
as $p_k(\mbox{\boldmath$\Omega $})
\in \Omega^{4k}(J^1\mathcal{M}_M)$.
Moreover, assuming $M$ is connected
and oriented, we have a principal $SO(n)$-bundle
over $J^1\mathcal{M}_M$, $O^{+}M
=\left\{
(u_{x},j_{x}^1g)\in OM\colon
u_{x}\,\text{is positively oriented}
\right\} $,
and $\mathbf{\omega }$ is
reducible to $O^{+}M$. A well-know result
(e.g., see \cite[Chapter 8]{GS},
\cite[XII, Theorem 2.7]{KN}) states that
$\mathcal{I}^{SO(n)}$ is generated by
the polynomials $\{p_k\}$ for odd $n$, and by
$\{p_k,\mathrm{Pf}\}$ for
even $n$, where $\mathrm{Pf}\in \mathcal{I}_{n/2}^{SO(n)}$
denotes the Pfaffian. For every even dimension $n=\dim M$,
we define the universal Euler form of $M$ by setting
\[
\mathrm{E}
=\left(
2\pi\right) ^{-\frac{n}{2}}
\mathrm{Pf}(\mbox{\boldmath$\Omega $})
\in \Omega ^n(J^1\mathcal{M}_M).
\]
From the identity
$\left(
2\pi
\right) ^{-n}\mathrm{Pf}^{2}=p_{n/2}$
we deduce
$\mathrm{E}\wedge\mathrm{E}
=p_{n/2}(\mbox{\boldmath$\Omega $})$.
The properties of $\mathbf{\omega }$
lead us readily to the following

\begin{proposition}
\label{pont}We have

\begin{enumerate}
\item [\emph{(1)}]The universal Pontryagin
forms and the universal Euler form are closed.

\item[\emph{(2)}] The universal Pontryagin forms
$p_k(\mbox{\boldmath$\Omega $})$
(resp.\ the universal Euler form $\mathrm{E}$)
are invariant under the action of
$\mathrm{Diff}M\times \mathbb{R}$ (resp.\
$\mathrm{Diff}^{+}M\times\mathbb{R}$) on
$J^1\mathcal{M}_M$.

\item[\emph{(3)}] For any Riemannian metric $g$
on $M$, we have
\begin{align*}
(j^1g)^\ast (p_k(\mbox{\boldmath$\Omega $}))
&
=p_k(\Omega^{g}),\\
(j^1g)^\ast (\mathrm{E})
&
=(2\pi )^{-\frac{n}{2}}\mathrm{Pf}(\Omega ^g),
\end{align*}
where $\Omega ^g$ denotes the curvature
of the Levi-Civita connection $\omega ^g$ of $g$.
\end{enumerate}
\end{proposition}

\begin{remark}
The Euler form is not invariant under the elements
of $\mathrm{Diff}^{-}M$ (recall that $O^{+}(M)$ is
invariant under $\mathrm{Diff}^{+}M$, but not under
$\mathrm{Diff}^{-}(M)$). In fact, for any
$\phi\in \mathrm{Diff}^{-}M$ we have
$\bar{\phi }^{(1)\ast }(\mathrm{E})=-\mathrm{E}$.
\end{remark}

The relation between the universal Pontryagin
and Euler forms on
$J^1\mathcal{M}_M$ and the usual Pontryagin
and Euler classes on $M$ is
the same as the relation between characteristic
forms and classes on the
bundle of connections of a principal bundle (e.g.,
see \cite{conn}): The map
$q_1^\ast \colon H^\bullet (M)\to H^\bullet
(J^1\mathcal{M}_M)$ is an isomorphism with inverse
map $(j^1g)^\ast \colon H^\bullet
(J^1\mathcal{M}_M)\to H^\bullet(M)$, for any Riemannian
metric $g$ on $M$, and by \ref{pont}--(3)
the $k$-Pontryagin (resp.\ Euler) class of $M$
is the image under this isomorphism
of the cohomology class of the universal $k$-Pontryagin
(resp.\ Euler) form.

As in the case of the bundle of connections
of a principal bundle, the Pontryagin forms contain
more information than the Pontryagin classes. For
example, if $4k>n$, the $k$-Pontryagin class vanishes,
but the corresponding form does not necessarily, as
$\dim J^1\mathcal{M}_M>n$. For example, if $n=2$, then
the first Pontryagin form $p_1(\mbox{\boldmath$\Omega $})
\in \Omega ^4(J^1\mathcal{M}_M)$ does not vanish,
whereas the first Pontryagin class vanishes
by dimensional reasons. According to \cite{equiconn},
these higher-order Pontryagin forms can be interpreted
as closed $\mathrm{Diff}^{+}M$-invariant differential
forms on the space of Riemannian metrics on $M$.
In a forthcoming paper, we shall study these forms
and their extension to equivariant cohomology
in a similar way as done in \cite{equiconn} for the
characteristic forms on the bundle of connections.

\section{The universal Levi-Civita connection characterized}

\begin{theorem}
\label{caracterizacion}
The universal Levi-Civita connection $\mathbf{\omega }$
is the only $\mathrm{Diff}M$-invariant connection form
on $q_1^\ast FM\to J^1\mathcal{M}_M$ satisfying the
condition \emph{\eqref{nablaG}}. In other words, the form
$\mathbf{\omega }$ is the only $\mathrm{Diff}M$-invariant
connection form on the bundle $OM\to J^1\mathcal{M}_M$.
\end{theorem}

The proof of this theorem is based in the following

\begin{lemma}
\label{lematecnico}
The $\mathrm{Diff}M$-invariant $1$-forms
on $J^1\mathcal{M}_M$ with values on
$\otimes ^2T^\ast M$ are $\lambda \theta
+\mu \mathrm{tr}\vartheta \otimes \mathbf{g}$,
$\lambda ,\mu \in \mathbb{R}$.
Equivalently, the only $\mathrm{Diff}M$-invariant
connection forms on the bundle
$q_1^\ast FM\to J^1\mathcal{M}_M$ are
$\mathbf{\omega }+\lambda \vartheta
+\mu \mathrm{tr}\vartheta \otimes \mathrm{id}_{TM}$,
$\lambda ,\mu \in\mathbb{R}$.
\end{lemma}

\begin{proof}
[Proof of Theorem \ref{caracterizacion}]
The universal Levi-Civita connection
$\mathbf{\omega }$ satisfies the conditions
of the statement by virtue of
Theorem \ref{propomeganegrita}. Conversely,
let us suppose that $\omega $
is another $\mathrm{Diff}M$-invariant connection
on $q_1^\ast FM$. Then, we have
$\omega =\mathbf{\omega }+\alpha $, with
$\alpha \in \Omega ^1
(J^1\mathcal{M}_M,\mathrm{End}TM)$.
Clearly, $\alpha $ is $\mathrm{Diff}M$-invariant,
and from Lemma \ref{lematecnico} we have
$\alpha =\lambda \vartheta
+\mu \mathrm{tr}\vartheta \otimes \mathrm{id}_{TM}$
for some $\lambda ,\mu \in \mathbb{R}$. Hence
$\nabla ^\omega \mathbf{g}=-2
\left(
\lambda \theta
+\mu \mathrm{tr}\vartheta \otimes \mathbf{g}
\right) $,
and consequently, $\omega $ satisfies the condition
\eqref{nablaG} if and only if $\lambda =\mu =0$;
i.e., if and only if $\omega =\mathbf{\omega }$.
\end{proof}

Next, we state some necessary results to prove
Lemma \ref{lematecnico}. First
of all, we recall the following:

\begin{theorem}
[\cite{ABP, Spivak}]\label{TFO(n)}
\emph{(Fundamental theorem of the invariant
theory for the orthogonal group) }
Let
$(v,w)\mapsto \left\langle v,w\right\rangle $
be the standard scalar product on
$V=\mathbb{R}^n$, allowing us to identify $V$
with its dual space. We consider the tensorial
representation of $O(n)$ on $\otimes ^kV$. Then,
we have

\begin{enumerate}
\item [\emph{(1)}] For $k$ odd, the unique
$O(n)$-invariant element of
$\otimes ^kV$ is the zero element.

\item[\emph{(2)}] For $k=2l$ even, the subspace
of $O(n)$-invariant elements
on $\otimes ^kV$ is generated by the following
invariant linear forms
\[
\varphi _{i_1,i_2,\dotsc,i_{2l-1},i_{2l}}
(v_1,\dotsc,v_{2l})
=\left\langle
v_{i_1},v_{i_2}
\right\rangle
\cdots
\left\langle
v_{i_{2l-1}},v_{i_{2l}}
\right\rangle ,
\]
where $i_1,i_2,\dotsc,i_{2l-1},i_{2l}$
stands for an arbitrary permutation
of the set of indices $1,2,\dotsc,2l-1,2l$.
\end{enumerate}
\end{theorem}

\begin{theorem}
[\cite{Spivak}]\label{TFSO(n)}
\emph{(Fundamental theorem of the invariant
theory for the special orthogonal group)}
With the previous notations, for
$k<n$, the $SO(n) $-invariants on
$\otimes ^kV$ coincide with $O(n)$-invariants.
For $k=n$, the space of $SO(n)$-invariants
is generated by the space of $O(n)$-invariants
and $\wedge ^nV$.
\end{theorem}

\begin{remark}
\label{contraccioncompleta}For $k=4$,
from Theorem \ref{TFO(n)} we conclude
that the $O(n)$-invariants are generated
by the following tensors:
\begin{align*}
\xi _1
&
=\sum \nolimits_{i,j}e_i
\otimes e_i
\otimes e_j
\otimes e_j,\\
\xi_{2}
&
=\sum \nolimits_{i,j}e_i
\otimes e_j
\otimes e_i
\otimes e_j,\\
\xi_{3}
&
=\sum\nolimits_{i,j}e_i
\otimes e_j
\otimes e_j
\otimes e_i.
\end{align*}
\end{remark}

\begin{proposition}
\label{inv}Let $(e_1,\dotsc,e_n)$
be the standard orthonormal base in
$V=\mathbb{R}^n$. Consider the $O(n)$-module
$E=\otimes ^3V\oplus
\left(
S^2V\otimes (\otimes ^2V)
\right)
\oplus
\left(
S^{2}V\otimes(\otimes ^{3}V)
\right) $.

\begin{enumerate}
\item [\emph{(1)}]The invariant elements
under the action of $O(n)$ on $E$ are
\begin{equation}
\eta
=\lambda \sum_{i,j}e_i
\odot e_i
\otimes e_j
\odot e_j
+\mu \sum _{i,j}e_i
\odot e_j
\otimes e_i
\odot e_j,
\quad
\lambda,
\mu \in \mathbb{R},
\label{contr}
\end{equation}
where $\odot$ denotes the symmetric product.

\item[\emph{(2)}] For $n\geq4$ the $SO(n)$-invariants
on $E$ coincide with the $O(n)$-invariants.
\end{enumerate}
\end{proposition}

\begin{proof}
(1) Every direct summand in $E$ is $O(n)$-invariant;
hence we need only to analyze the invariants on each
summand.

From Theorem \ref{TFO(n)} it follows that there are
no invariants on $\otimes ^3V$ and on
$S^2V\otimes (\otimes ^3V)$, and the invariants on
$S^2V\otimes (\otimes ^2V)$ are obtained by linear
combination of the elements cited in
Remark \ref{contraccioncompleta}.
Hence, they are of the form
$\xi =\lambda \xi _1+\mu \xi _2+\nu \xi _3$.
Then,
$\xi \in S^2V\otimes (\otimes ^2V)$
if and only if $\mu =v$, and we obtain \eqref{contr}.

(2) For $n>5$ the result follows from Theorem \ref{TFSO(n)}.
For $n=5$, from Theorem \ref{TFSO(n)} it follows that
$O(n)$-invariants and $SO(n)$-invariants coincide on
$\otimes ^3V$ and on $S^2V\otimes (\otimes ^2V)$.
Moreover, since we have
$\wedge ^5V\cap (S^2V\otimes (\otimes ^3V))=0$,
the same conclusion holds for the remaining summand
$S^2V\otimes (\otimes ^3V)$. Finally, for $n=4$,
again from Theorem \ref{TFSO(n)} and the fact that
$\wedge ^4V\cap(S^2V\otimes (\otimes ^2V))=0$,
it follows that there are no new invariant on
$\otimes ^3V$ or on $S^2V\otimes (\otimes ^3V)$. Also,
as $-\mathrm{id}_V\in SO(4)$ and
$(-\mathrm{id}_V)\cdot(\eta )=-\eta $ for all
$\eta \in S^2V\otimes (\otimes ^3V)$, we conclude that
no new invariant appears on $S^2V\otimes (\otimes ^3V)$,
and the result follows.
\end{proof}

\begin{proof}
[Proof of Lemma \ref{lematecnico}]
Let us fix a point
$z_0=j_{x_0}^1g_0\in J^1\mathcal{M}_M$,
and let us consider a normal coordinate system
$(U;x^i)$ centered at $x_0$ for the metric $g_0$.
The expression of a covector
$\eta \in \Omega ^1
(J^1\mathcal{M}_M,\otimes ^2T^\ast M)$
at $z_0$ on this coordinate system is
$\eta _{z_0}
=(\lambda _{ab,i}dx^i
+\lambda _{ab}^{ij}dy_{ij}
+\lambda _{ab}^{ijk}dy_{ij,k})_{z_0}
\otimes (dx^a)_{z_0}\otimes(dx^b)_{z_0}$.

For reversible $M$ we set $G=SO(n)$ and for irreversible
$M$, $G=O(n)$. Given $A\in G$, we define a local
diffeomorphism $\varphi _A\colon U\to M$
around $x_0$ by $\varphi _A^i(x)=A_j^ix^j$,
$\forall x\in U$. As $\varphi _A$ is a linear
transformation, from the expression \eqref{coor}
we deduce $\overline{\varphi _A}(j_{x_0}^1g_0)
=j_{x_0}^1g_0$ and we have
\begin{align*}
\overline{\varphi _A}^{(1)\ast }(dx^i)_{z_0}
&
=(A^{-1})_a^i(dx^a)_{z_0}
=\sum \nolimits _aA_i^a(dx^a)_{z_0},\\
\overline{\varphi _A}^{(1)\ast }(dy_{ij})_{z_0}
&
=A_i^aA_j^b(dy_{ab})_{z_0},\\
\overline{\varphi _A}^{(1)\ast }(dy_{ij,k})_{z_0}
& =A_i^aA_j^bA_k^{c}(dy_{ab,c})_{z_0}.
\end{align*}
Hence the map
\begin{align*}
(dx^i)_{z_0}
&
\longmapsto e_i,\\
(dy_{ij})_{z_0}
&
\longmapsto e_i
\odot e_j,\\
(dy_{ij,k})_{z_0}
&
\longmapsto e_i
\odot e_j
\otimes e_k,
\end{align*}
determines a $G$-module isomorphism between
$T_{z_0}^\ast J^1\mathcal{M}_M
\otimes (\otimes ^2T_{x_0}^\ast M)$
and the space $E$ in the statement
of Proposition \ref{inv}. The local diffeomorphism
$\varphi _A$ satisfies the conditions in
Lemma \ref{dext} and, hence, there exists
$\phi _A\in \mathrm{Diff}M$ extending
$\varphi _A$ on a neighbourhood of $x_0$.
As $\eta $ is $\mathrm{Diff}M$-invariant we have
$\overline{\phi _A}^{(1)\ast}(\eta )=\eta $,
and hence $\overline{\varphi _A}^{(1)\ast }
\left(
\eta _{z_0}
\right)
=\eta _{z_0}$ for every $A\in G$.
As $n=\dim M\geq 4$ for
an irreversible $M$, from Proposition \ref{inv},
for some $\lambda ,\mu \in \mathbb{R}$, we obtain
\begin{align*}
\eta _{z_0}
\!\!
&
=\!\!\lambda \sum \nolimits _{i,j}(dy_{ij})_{z_0}
\otimes (dx^i)_{z_0}\otimes (dx^j)_{z_0}
+\mu \sum \nolimits _{i,j}(dy_{ii})_{z_0}
\otimes (dx^j)_{z_0}\otimes (dx^j)_{z_0}\\
\!\!
&
=\!\!
\lambda \theta _{z_0}
+\mu
(\mathrm{tr}\vartheta \otimes \mathbf{g})_{z_0},
\end{align*}
As the point $z_0\in J^1\mathcal{M}_M$
is arbitrary, for certain smooth functions
$a,b$ on $J^1\mathcal{M}_M$ we have
$\eta =a\theta
+b\mathrm{tr}\vartheta \otimes \mathbf{g}$.
As $\eta $, $\theta $ and
$\mathrm{tr}\vartheta \otimes \mathbf{g}$
are $\mathrm{Diff}M$-invariant, $a$ and $b$ are
also $\mathrm{Diff}M$-invariant and, by virtue
of Proposition \ref{trans}, they are constant.
\end{proof}

\begin{remark}
The characterization of the connection
$\mathbf{\omega }$ given on Theorem \ref{caracterizacion}
does not hold for higher-order jet bundles. In fact,
below we sketch the proof of the existence of a natural
$1$-form $\alpha \in \Omega ^1
(J^3\mathcal{M}_M,\mathrm{End}_{\text{\textsc{A}}}TM)$.
Hence, $\mathbf{\omega }+\alpha $ is a
$\mathrm{Diff}M$-invariant connection form on
$q_3^\ast FM$, which also satisfies the condition
\eqref{nablaG}, by virtue of Lemma \ref{lemasym}. Let
\begin{multline*}
\theta ^3
=\left(
dy_{ij}-y_{ij,r}dx^r
\right)
\otimes \partial /\partial y_{ij}
+\left(
dy_{ij,k}-y_{ij,kr}dx^r
\right)
\otimes \partial /\partial y_{ij,k}\\
+\left(
dy_{ij,kl}-y_{ij,klr}dx^r
\right)
\otimes \partial /\partial y_{ij,kl},
\quad i\leq j,k\leq l,
\end{multline*}
be the $(q_{32})^\ast V(q_{2})$-valued
$1$-form on $J^3\mathcal{M}_M$
defining its contact structure (cf.\ \cite{Mu}).
We first notice the natural exact sequence
of vector bundles over $J^2\mathcal{M}_M$,
\[
0\!\to \!
\left(
q_{2}
\right) ^\ast
\left(
S^2T^\ast M\otimes
S^2T^\ast M
\right)
\!
\overset{\iota _2}{\longrightarrow }
\!
\left(
q_2
\right) ^\ast J^2
\left(
S^2T^\ast M
\right)
\!
\to
\!
\left(
q_2
\right) ^\ast
J^1
\left(
S^2T^\ast M
\right)
\!
\to
\! 0
\]
splits naturally, as a retract
$\rho _2\colon
\left(
q_2
\right) ^\ast J^2
\left(
S^2T^\ast M
\right)
\to
\left(
q_2
\right) ^\ast
\left(
S^2T^\ast M\otimes S^2T^\ast M
\right) $
exists of $\iota _2$
given by,
\begin{multline*}
\rho _2
\left(
j_x^2g,j_x^2h
\right)
(X_1,X_2,X_3,X_4)
=\tfrac{1}{2}
\left(
\nabla ^g
\right) _x^2(h)
(X_3,X_4,X_1,X_2)\\
+\tfrac{1}{2}
\left(
\nabla ^g
\right) _x^2(h)
(X_4,X_3,X_1,X_2),
\end{multline*}
for all
$j_x^2g\in J_x^2
\left(
\mathcal{M}_M
\right) $,
$j_x^2h\in J_x^2
\left(
S^2T^\ast M
\right) $,
and $X_1,\dotsc,X_4\in T_xM$.
Let $c_{24}\colon
\left(
q_{2}
\right) ^\ast
\otimes ^4T^\ast M\to
\left(
q_{2}
\right) ^\ast
\otimes ^2T^\ast M$
be the metric contraction of the second
and fourth arguments, i.e.,
$c_{24}(j_x^2g,X_1
\otimes X_2
\otimes X_3
\otimes X_4)
=(j_x^2g,g(X_2,X_4)X_1
\otimes X_3)$.
By using the canonical vector-bundle
isomorphism $V(q_2)\cong
\left(
q_2
\right) ^\ast
J^2
\left(
S^2T^\ast M
\right) $,
the form we are looking for, is defined as follows:
$\alpha =(c_{24}\circ \rho _2
\circ \theta ^3)_{\text{\textsc{A}}}$.
\end{remark}

\end{document}